\documentclass[12pt]{article}

\usepackage{amsmath,amssymb, amscd}

\newcommand {\BOX} {\rule{2mm}{2mm}}

\newcommand {\be} {\begin{equation}}
\newcommand {\eeqn} {\end{equation}}
\newcommand {\bea} { \begin{eqnarray}}
\newcommand {\eea} {\end{eqnarray}}
\newcommand {\beas} { \begin{eqnarray*}}
\newcommand {\eeas} {\end{eqnarray*}}
\newcommand {\ra} {\rightarrow}

\newtheorem {lemma} {LEMMA} [section]
\newtheorem {thm}[lemma]{THEOREM}
\newtheorem {prop}[lemma]{PROPOSITION}

\newtheorem {cor}[lemma]{COROLLARY}

\newtheorem {conj}[lemma]{CONJECTURE}

\numberwithin{equation}{section}

\begin{document}


\title{On surfaces of general type
with maximal Albanese dimension\\
\large{\it }}
\author{Steven S. Y. Lu
\thanks{Partially supported by an NSERC grant and UQAM}}
\date{}

\maketitle

\begin{abstract}
Given a minimal surface $S$ equipped with a generically finite map to an Abelian
variety and $C\subset S$ a rational or an elliptic curve, we show that the 
canonical degree of $C$ is bounded by four times the self-intersection of
the canonical divisor of $S$. As a
corollary, we obtain the finiteness of rational and elliptic curves
with an optimal uniform bound on their canonical degrees on 
any surface of general type with two linearly independent 
regular one forms.
\end{abstract}

\section{Introduction}

The object of this paper is to give an effective bound on the canonical 
degrees of rational and elliptic curves on a minimal surface of 
general type with generically finite 
Albanese map linearly in terms of the canonical volume of the surface. 
Here, the canonical degree is the degree with respect to 
the canonical polarization of the surface and
the self-intersection of this polarization is the canonical volume.

We will work in the complex analytic setting so that all varieties
are complex analytic. We will assume the rudiments from the the theory of 
classification of complex projective surfaces as found for example
in \cite{BPV}. We call a smooth projective surface to be of 
maximal Albanese dimension if its Albanese
map is generically finite (N.B. This terminology is different from 
that introduced in \cite{Ueno}).    We recall that
any smooth projective surface $S$ of maximal Albanese dimension
admits a morphism to a minimal one whose
canonical divisor pulls back to the positive part $P$ of the
Zariski decomposition of the canonical divisor $K$ of $S$ and
the canonical volume of $S$ is 
by definition $vol(K):=P^2$.  Our main theorem is as follows.

\begin{thm}\label{main}
Let $S$ be a smooth projective surface with maximal Albanese 
dimension and $C$ a rational or an elliptic curve in $S$.
Let $P$ be the positive part in the Zariski decomposition 
of the canonical divisor $K$ of $S$. 
Then $$PC\leq 4\, vol(K).$$
\end{thm}

We remark that Y. Miyaoka has recently
given an effective canonical bound for surfaces with positive
Segre class \cite{Miya07} generalizing the same given by him and 
the author for surfaces with positive topological index \cite{LM1}.
It is easy to construct examples of surfaces with maximal
Albanese dimension which do not have positive Segre class. Prior to the present 
article, all such canonical bounds were obtained via the powerful but 
highly nontrivial log-orbifold Miyaoka-Yau inequality on surfaces, 
for which
we cite \cite{Bog, Miya, Yau, Sakai, LM1, LM2, Langer, Miya07}.
But the bounds so obtained are neither simple
nor are optimal. In this article,
we introduce an elementary approach that
relates geometrically the canonical degree of curves 
directly with the canonical class of the surface.
In view of the simplicity and aptitude
of the idea, the method introduced
and the resulting bound,
we expect it to have significance and applicability
to arbitrary algebraic curves and 
even to holomorphic curves as well as to arithmetic objects
as it seems to relate to the
standard conjectures in these subjects; its 
generalizations to arbitrary curves 
would have strong
implications, for example, for the abc-conjecture over 
function field for an abelian variety.

We also remark that a weakened version of our theorem has already
been given by Noguchi-Winkelmann-Yamanoi \cite{NWY07} where the Albanese
map is assumed to be finite and surjective, 
that is they exclude the possibility of
exceptional fibers in their hypothesis although
they do deal with the quasi-projective case at
the same time.  In that case, they 
bound the canonical degree of elliptic curves in terms of the degree of 
Gauss map of the canonical divisor of the surface, which they
do not bound. Our bound 
gives a sharp effective bound for the degree of this Gauss map
and we expect it to be sharp also for the canonical degree of
elliptic curves on such a surface, if not extremely close to it,
being by its very nature the best bound 
allowed by this approach.

We finally remark that it is possible to obtain less effective 
(noneffective) bounds more simply. But as all the details presented
here are quite general in nature and has much wider applicability,
we believe that they are worth the trouble for the optimal
bound given.\\

Since a surface of general type cannot support a nontrivial family of
rational or elliptic curves and a bound on the canonical degree of 
such curves puts them in a bounded family, 
our result implies that that there are only
a finite number of rational and elliptic curves on a surface of 
general type with generically finite Albanese map. In fact, our 
result is strictly stronger as it implies a global bound on such
finiteness, in a smooth family of such surfaces for example.  
More generally, for a surface
with irregularity at least two, if its Albanese map is not generically
finite, then it admits a map to a hyperbolic curve by the structure 
theorem of Kawamata and Ueno \cite{Kaw}. Hence rational and
elliptic curves lie on the fibers of such a map and are thus
finite in number if the surface is of general type. It follows
that there are only a finite number of rational and elliptic 
curves on any surface of general type with irregularity at
least two. The following is an immediate 
corollary of this
and of the main theorem of Noguchi-Winkelmann-Yamanoi in
\cite{NWY} giving the algebraic degeneracy of holomorphic
curves in the remaining case of such surfaces.

\begin{cor}
Let $S$ be a smooth projective surface of general type with irregularity
two or more. Then $S$ admits a 
proper Zariski subset that contains all 
nontrivial holomorphic images of $\Bbb C$.
\end{cor}

In general we have the following sweeping conjecture 
concerning the algebraic pseudo-hyperbolicity of varieties of
general type. The conjecture is at
least well indicated from the works of F. Bogomolov, B. Mazur,
M. Green and P. Griffiths in the seventies and that of Y. Kawamata in the 
early eighties and made explicit by Serge Lang in 1986 \cite{Lang}:

\begin{conj} An algebraic variety of general type admits 
a proper Zariski subset that contains all subvarieties 
not of general type.
\end{conj}

Aside from Bogomolov's result validating the conjecture for
surfaces of positive Segre class \cite{Bog}, the best general evidence
for such a conjecture up till now was given by Kawamata \cite{Kaw1}
validating the conjecture for subvarieties of Abelian varieties in
characteristic zero. This result was later generalized to the case
of semi-Abelian varieties by Noguchi \cite{Nog86} and to the case
when the field of definition is an 
algebraically closed field of arbitrary characteristic \cite{Ab}.
However, none of these result is effective and, besides the above 
mentioned result of Noguchi-Winkelmann-Yamanoi, the conjecture 
seemed unknown in general for a variety with maximal Albanese dimension.
Our method here applies directly to this latter case for its affirmative
resolution \cite{Lu10}.\\

A short but already quite telling description of the proof of 
this main theorem 
is as follows. The theorem reduces easily to the case when
$S$ is minimal (so that $P=K$) and 
admits a surjective morphism $\alpha$ to an abelian 
surface $A$
and to the case $C$ is neither $\alpha$-exceptional,
so that $C_0=\alpha(C)$ is an elliptic curve in $A$, 
nor contained in the
ramification locus of $\alpha$. In this case, 
$(\det d\alpha)$
is a canonical choice for $K$ and, outside its
common divisorial locus with $\hat C=\alpha^*(C_0)$,
$K$ intersects properly with 
$\hat C$ and the relevant local intersection numbers with $C$
are dominated by the local intersection numbers of $\hat C$
with each such component 
of $K$.  These local intersection numbers along a
horizontal component of 
$D=K_{red}$ can be interpreted in terms of the 
vanishing degree of the natural section of $K+D$ along $D$ 
given by $(\alpha^*w_0)\wedge ds/s$ where $w_0$ is the nonzero
holomorphic one form on $A$ such that 
$TC_0\in \ker w_0$. Since the intersection
number of $\hat C$ with any vertical component of $D$ 
or with any horizontal component of $D$ in $\hat C$ is zero, 
finding a decomposition of $K$ into parts having
nonnegative intersection number with $K+D$ 
and summing yield the result.

As for the plan of the paper,
section two deals with this decomposition of $K$
while section three gives the key lemmas
for the proof of the main theorem in 
section 4.\\

The author is greatful to B. Shiffman, 
J. Noguchi, J. Winkelmann and A. Fujiki 
for their interests in this work and 
for their kind invitation to Japan where this work
was presented.

\section{A decomposition of divisors on surfaces}

We will let $S$ be a smooth projective surface and $K$ a
canonical divisor for $S$ throughout. Let $D, D'$
be two divisors in $S$. We will denote $D\leq D'$ to mean that $D'-D$ is
effective (possibly zero). By a component of a divisor, we will mean a
prime component as opposed to a connected component of a divisor, so 
that it is reduced and irreducible.

\begin{lemma}Let $D$ be a reduced divisor in $S$ and $C$ a connected
divisor contained in $D$, i.e. $C\leq D$. If $(K+D)C<0$, then either 
$(K+D)C=-2$ or $(K+D)C=-1$. These two cases, which we will call
case 1 and case 2, imply respectively: 
\begin{enumerate}
\item $(K+C)C=-2$ and $(D-C)C=0$. The latter means that $C$ is isolated 
in $D$, i.e. $C$ is not connected to the rest of $D$.
\item $(K+C)C=-2$ and $(D-C)C=1$. The latter means that $C$ intersect
the rest of $D$ at one point (and transversely). 
\end{enumerate}
\end{lemma}

\noindent{\bf Proof:} Since $C$ is connected, its arithmetic genus $g$ is
given by $1-(h^0({\cal O}_C) - h^1({\cal O}_C))=h^1({\cal O}_C)\geq 0.$
Hence by adjunction, $(K+C)C=2g-2$ is an even number that is at least
$-2$. The rest follows from $(K+D)C=(K+C)C+(D-C)C$. \BOX\\

Recall that a tree of curves is a connected reduced curve $C$ which
is disconnected by removing any of its component curves. Any component
of such a tree is connected to any other component by a unique path.
A rooted tree of curves is a tree with a component curve singled out, 
called the $root$. Such a tree becomes directed by the partial ordering 
given by $C_i\leq C_j$ if the unique path from the root to the
component $C_j$ passes through $C_i$, in which case we also say
that $C_j$ is a descendant of $C_i$. Note in particular that every
component curve is a descendant of itself in a directed tree.
A component curve $C_k$ without a
descendant beside itself
is called a $leave$. In a rooted tree of curves, the union
of a component curve with all
its descendant is call a $branch$. A tree without a root specified 
is called a free tree and such a tree becomes directed by 
either choosing a root curve or by choosing the 
unique transversal
intersection of a pair of curves if exists, 
with the obvious ordering given
as before. A disjoint collection of trees is called
a forest.

\begin{lemma}\label{L1}
In the previous lemma, a curve $C$ in case 1 consists of an isolated (free) 
tree of rational curves intersecting transversally 
and a curve $C$ in case 2 consists of a rooted tree of rational curves 
intersecting transversally. If we direct the tree of case 1 by choosing a 
transversal intersection point of two 
curves if it exists, then every branch of the tree
in either case satisfies the conditions of case 2.
\end{lemma}

\noindent{\bf Proof:} We will only show the lemma in case 2 since case 1
follows similarly. In this case, $(K+D)C=-1$ and $C$ intersects the rest of
$D$ at one point $p$. Let $C_0$ be the component of $C$ containing $p$ and
let $C_{01},..., C_{0k}$ be the connected components of $C-C_0$. Let
$D'=D-\sum_i C_{0i}$.
By construction, $(K+D)C_{0i}\geq -1$ and $C_{0i}(C-C_{0i})\geq 1$ for all $i$
and since $C_0$ is a connected nonisolated component of $D'$, we have 
by the previous lemma that
\begin{eqnarray*}
-1&\leq& (K+D')C_0= (K+D)C_0-\sum_i C_{0i}C_0\\&=&
(K+D)C-\sum_i[(K+D)C_{0i}+C_{0i}(C-C_{0i})]\leq (K+D)C=-1. 
\end{eqnarray*}
This forces equality everywhere and so $(K+D)C_{0i}= -1$ for all $i$
and $(K+D')C_0=-1$. This means that the $C_{0i}$'s are all 
curves belonging to case 1 connected to $C_0$ and that $C_0$ 
is a smooth rational curve as $(K+C_0)C_0=-2$ by the previous lemma.
Hence we see that $C$ is a rooted tree of rational curves by induction, 
with $C_0$ as a root. \BOX\\

 We now set
$$F(D)=\{\ C\leq D\ |\ C\ \mbox{is connected and}\ (K+D)C<0\ \}.$$ 
It follows that the union of the elements in $F(D)$ form 
a subset of $D$ that is a forest. We will call this
forest $F_D$. It is a disjoint union of maximal elements in 
$F(D)$ and the leaves in this forest
are exactly the prime components of $D$ belonging to $F(D)$.\\

We apply our result to the case $D=K_{red}$, the reduced canonical
divisor, to obtain a positive partial decomposition of $K$. 
The following
lemma is the key to this decomposition.\\

Recall that a curve
$C$ is called  a Hirzebruch-Jung string 
if it is a rooted tree of curves with exactly one
leave so that we have actually a total ordering on $C$.
In this case, we may write the prime decomposition of 
$C$ as $C_0+C_1+\cdots + C_r$
where $C_iC_j$ is one or zero according to whether $j=i+1$
or not  and $C_0$, $C_r$ are the root and leave of the string
respectively. We will denote the subset of $F(D)$ consisting
of elements that are not isolated in $D$ and that are
Hirzebruch-Jung strings by $F^{HJ}(D)$. 

\begin{lemma}\label{L2} Let notation be as above but with $D=K_{red}$.
Let $C\in F^{HJ}(D)$  
with the prime
decomposition $C_0+C_1+\cdots + C_r$ as given above. Let $C_{-1}$
be the unique component of $D-C$ connected to $C$. Let
$m_i>0$ be the multiplicity of $K$ on $C_i$. If $S$ is minimal,
then for $i=-1,0,\dots,r-1$, 
we have 
$$m_i\geq 2m_r.$$
\end{lemma}

\noindent{\bf Proof:} Since each component $C_i$ of $C$ is a 
rational curve and $S$ is minimal, we have, for each $i\geq 0$, 
$$C_i^2=-2-KC_i\leq -2.$$
Let $\hat C=C_{-1}+C$. Then 
$$(m_{r-1}+1)+(m_r+1)C_r^2=(K+D)C_r=-1,$$
and we obtain $m_{r-1}+2=-(m_r+1)C_r^2\geq 2(m_r+1)$. 
Hence $m_{r-1}\geq 2m_r$ and our lemma is true for $i=r-1$.
We now establish the lemma by showing that $m_{i-1}>m_i$ for
all $i\geq 0$. This being already established for $i=r$, we
assume by induction that $m_i>m_{i+1}$ for an $i$ strictly between
$-1$ and $r$. For such an $i$, we have by lemma~\ref{L1} that
$$0=(K+D)C_i=(m_{i-1}+1)+(m_i+1)C_i^2+ (m_{i+1}+1)< m_{i-1}+1 +(m_i+1)(C_i^2+1),$$
and so $m_{i-1}>m_i$ and the result follows by induction. \BOX\\

\noindent
We remark that if $S$ is not minimal, the lemma would still hold
with a weaker inequality.\\

\begin{lemma}\label{L3}
With the notations as before and $D=K_{red}$, let $C\in F(D)$
be a maximal element non-isolated in $D$. Let $C=\sum_{i\in I} C_i$ be the 
prime decomposition of $C$ and let $K_C=\sum_{i\in I} m_iC_i$ where
$m_i>0$ is the multiplicity of $K$ on $C_i$. If $S$ is minimal
and $C$ is not a Hirzebruch-Jung string, then $(K+D)K_C\geq 0.$
\end{lemma}

\noindent{\bf Proof:} We first introduce some terminology. Let
$C\in F(D)$ be non-isolated in $D$. Let $C_i$ be a component of $C$.
Recall that $\bar C_i$ is the union of all the descendant of $C_i$.
A component $C_i$ of $C$ is called split 
if $\bar C_i -C_i$ contains more than one connected components, 
at least one of which
lies in $F^{HJ}(D)$. In this case,
we let $K_i= m_iC_i+\sum_j m_{ij}C_{ij}$ where the $C_{ij}$'s are
the connected components of $\bar C_i -C_i$ lying in $F^{HJ}(D)$,
$m_i$ is the multiplicity of $K$ on $C_i$ and $m_{ij}$ is the 
multiplicity of $K$ on the leave of $C_{ij}$. Since 
$(K+D)\bar C_i=-1$, $(K+D)C_{ij}=-1$ for all $j$ and  
$(K+D)(\bar C_i-C_{i1})\leq -1$ by Lemma~\ref{L1} in this case, 
we obtain easily from Lemma~\ref{L2} that
$$(K+D)K_i\geq 0.$$
By assumption, $C$ is not a Hirzebruch-Jung string. So every leave
of $C$ is the descendant of a split component of $C$. 
Let $C_i,\ i\in I'$ be the collection of split components of $C$
and $K_i$ be as above. Now $K'=\sum_{i\in I'}K_i\leq K_C$ by lemma~\ref{L2}
and $K''=K_C-K'$ has zero multiplicity on the leaves by construction. 
Hence $K''$ has
only support on the $C_j$'s which are not leaves and as
$(K+D)C_j\geq 0$ for such $C_j$'s, $(K+D)K''\geq 0$. It follows that
$$(K+D)K_C=(K+D)K''+\sum_{i\in I'}(K+D)K_i\geq 0.\ \BOX$$

\section{The key proposition}

Let $G$ be a reduced 
irreducible divisor on a smooth projective
surface $S$ 
and defined by the section
$s\in H^0(S, {\cal O}(G))$, i.e., $(s)=G$. 
The key point to our approach is that
$ds/s$ gives a well defined holomorphic section 
of $\Omega_S(G)$ over $G$. 
We will only need this fact in the following situation.

\begin{lemma} With the data as given above, let $w$ be a
holomorphic one-form on $S$.
Then $w\wedge ds/s$ gives a well defined 
holomorphic section of ${\cal O}(K+G)|_G$, i.e.,
$$u=w\wedge ds\in H^0(G,{\cal O}(K+G)).$$
If $w$ does not pull back to the zero one-form on any component of $G$,
then $u$ is nowhere identically vanishing.
\end{lemma}

\noindent{\bf Proof:} Given two local trivializations for 
${\cal O}_G(G)$ on $G\cap U$ where $U$ is
a Stein neighborhood of a point $p\in G$, consider
two  trivializations of ${\cal O}_U(G)$ that extend them.
Let $s_1$ and $s_2$ be the respective holomorphic functions on 
$U$ representing $s$ with respect to the trivializations.
Then $t=s_1/s_2$ is a nonvanishing holomorphic function on 
$U$. Since  $ds_1=tds_2+s_2dt$ on $U$, we have 
$w\wedge ds_1=tw\wedge ds_2$ on $G\cap U$. It follows that
$w\wedge ds$ transforms as a holomorphic
section of $\Omega_S^2(G)={\cal O}(K+G)$
over $G$ under different local trivializations of ${\cal O}_G(G)$. 
The last part of the lemma is clear. \BOX\\

An important remark for the application in general is that
if $s$ is only required to satisfy $G=(s)_{red}$, then $ds/s$
would give a nowhere-vanishing section of $\Omega_S(G)$ over 
$G$ thereby effectively reducing the multiplicities of $(s)$
to one.

We now apply this lemma to obtain the key proposition
for the proof of our theorem.

\begin{prop} \label{Pr1}
Let $S$ be a smooth non-Abelian surface, 
$A$ a complex  Abelian surface and $\alpha: S\rightarrow A$ a 
surjective morphism. We consider the following data on $S$:
\begin{itemize}
\item[-] $K=(\det d\alpha)$ the canonical divisor of $S$ determined
by $\alpha$
\item[-] $D=K_{red}$ the reduced canonical divisor
\item[-] $G$ a horizontal component of $D$, i.e., $G$ is not 
$\alpha$-exceptional 
\item[-] $C_0$ an elliptic curve in $A$ considered as a reduced
divisor
\item[-] $\hat C=\alpha^*(C_0)$ the total transform of $C_0$
\item[-] $C'$ the sum of non-elliptic
horizontal components of $\hat C$ on which $\hat C$ has
multiplicity one
\item[-] $\bar C=\hat C - C' - C''$ 
where $C''$ is some vertical part of $\hat C$, i.e., 
$C''\leq \hat C$ is $\alpha$-exceptional

\end{itemize}
Assume that $G$ is not a component of $\bar C$. 
Then we have $$ (G,\bar C)\leq 2 (G, P+D) -n_G,$$
where $n_G$ is the intersection number of $D-G$ with
the smooth part of $G$.
\end{prop}

\noindent{\bf Proof:} Our assumption on $S$ implies easily that $\alpha$ 
must be ramified somewhere and so $\det d\alpha$ gives a well defined
divisor on $S$. We note also that, since $\hat C$ cannot be reduced on 
any vertical component of the ramification divisor,  our assumption on $G$
implies that it is not a component of $\hat C$ so that $G\cap (\hat C)_{red}$ 
is a finite set containing $Q=G\cap (\bar C)_{red}$ . 
We first assume that $G$ is smooth for the proof.\\

Let $w_0$ be the constant holomorphic 
one-form that defines the direction of $C_0$, i.e., locally
we can write $w_0=df_0$ with $(f_0)=C_0$. By definition,
$$ (G,\bar C)=\sum_{x\in Q} (G, \bar C)_x \leq
\sum_{x\in Q} (G, \hat C)_x.$$ 
Let $x\in Q$ and $U$ a Stein neighborhood of $\alpha(x)$.
Then $w_0|_U=df_0$ for a holomorphic function $f_0$ on $U$
with $(f_0)=C_0|_U$. Let $V=\alpha^{-1}(U)$ and 
$f=\alpha^*f_0=f_0\circ\alpha\in {\cal O}_V$. 
Then $(G, \hat C)_x$ is  by definition the vanishing order 
at $x$ of $f|_G$, i.e.,
$$(G, \hat C)_x=\mbox{ord}_x (f|_G).$$
Let $w=\alpha^*w_0$ and $\gamma_1=w\wedge ds|_G\in H^0(G, K+G)$. 
Then $w|_V=df$ and so
$$\mbox{ord}_x(\gamma_1)=\mbox{ord}_x(df|_G)=(G, \hat C)_x-1.$$
Let $G'$ be the reduced divisor $D-G$ considered also as a 
subset of $S$. Let $r$ be a section of ${\cal O}(G')$ with
$(r)=G'$ and let $\gamma=\gamma_1\otimes r$. Then $\gamma$
is a section of ${\cal O}(K+D)$ over $G$ and
\begin{equation}\label{ineq2}
\mbox{ord}_x(\gamma)=\mbox{ord}_x(\gamma_1)+(G, G')_x\geq (G,D-G)_x.
\end{equation}
\noindent
We observe that if $x\not\in G'$, then 
$(G,G')_x=0$ and $x$ lies in some horizontal 
component(s) of $\bar C$ away from any intersection point
with vertical 
components of $\hat C$. In this case, either $\bar C$ is not 
reduced at $x$, in which case 
$(G, \hat C)_x\geq (G, \bar C)_x\geq 2$,
or the components of $\bar C$ through $x$ 
are horizontal elliptic curves. In the latter case,
since the tangent directions of these elliptic curves 
cannot lie in the kernel of $d\alpha$ as no elliptic
curve can ramify over an elliptic curve, $d\alpha$ has
rank one at $x$ and maps the tangent "directions" of $C$ at $x$, 
and therefore also $T_x  S$, to 
$T_{\alpha(x)} C_0=\ker w_0(\alpha(x))$. Hence, in this
case, $w=\alpha^*w_0$ is identically zero on $T_x S$
and so $\mbox{ord}_x(\gamma_1)\geq 1$.
It follows that if $x\not\in G'$, then 
$\mbox{ord}_x(\gamma)=\mbox{ord}_x(\gamma_1)\geq 1$ and so
$$(G, \hat C)_x=\mbox{ord}_x(\gamma_1)+1\leq 2\,\mbox{ord}_x(\gamma_1)
=2\,\mbox{ord}_x(\gamma)-(G,D-G)_x.$$
On the other hand, if $x\in G'$, then $1\leq (G,G')_x$, and so we have again
$$(G, \hat C)_x=\mbox{ord}_x(\gamma_1)+1\leq \mbox{ord}_x(\gamma)
\leq 2\,\mbox{ord}_x(\gamma)-(G,D-G)_x.$$

\noindent Since the inequality~(\ref{ineq2}) is actually true for
all $x\in G$ by the construction of $\gamma$, we have 
$2\,\mbox{ord}_x(\gamma)-(G,D-G)_x\geq 0$ for all $x\in G$.
It follows that
$$(G, \bar C)\leq \sum_{x\in G}\big( 2\,\mbox{ord}_x(\gamma)-(G,D-G)_x\big)=
2 (G, K+D) - (G, D-G).$$
This proves the proposition in the case $G$ is smooth.\\

In the case $G$ is not smooth, let $\pi:\tilde S\ra S$ be a minimal
resolution of $G$. By replacing $G$ with its strict transform $\tilde G$, 
$\bar C$ with its total transform $\pi^*\bar C$, $K$ with the canonical
divisor $\tilde K=(\det d (\alpha\circ \pi))$ of $\tilde S$ and 
$D$ with $(\tilde K)_{red}=(\pi^*D)_{red}$ we find that all the assumptions
of the proposition still hold and therefore
$$(G, \bar C)=(\tilde G, \pi^*\bar C)\leq 
2 (\tilde G, \tilde K +\tilde D) - (\tilde G, \tilde D - \tilde G)
\leq 2 (\tilde G, \tilde K +\tilde D) - n_G,$$
the last inequality owing to the fact that 
the intersection number $n_G$ of $D-G$ with the smooth part 
of $G$ is unaffected by $\pi$ and so can be no greater than 
$(\tilde G, \tilde D - \tilde G)$. The proposition now follows from
the claim that $(\tilde G, \tilde K + \tilde D)\leq (G,K+D)$
for which it is sufficient via induction to verify for the case
$\pi$ is a single blow-up. But this case follows directly from 
$\tilde K=\pi^*(K)+E$ and $\tilde D=\pi^*(D)-(m-1)E$ where
the multiplicity $m$ of the point of $D$ blown up is 
necessarily no less than two as $D$ is singular there. \BOX

\section{Proof of the main theorem}

We first reduce the main theorem, Theorem~\ref{main}, to the case when the 
surface $S$ is minimal of general type admitting a surjective morphism to
an Abelian surface. \\

Recall that the Albanese map of a smooth projective variety 
$X$ is a morphism $\alpha_X: X\ra A$
where $A$ is an Abelian variety
(i.e., $A$ is a projective variety whose universal covering is ${\Bbb C}^n$ 
for some $n>0$)
and where the pair $(\alpha_X, A)$ is characterized by the following universal 
property: Any morphism $\beta : X\ra A'$ where $A'$ is an Abelian
variety admits a unique factorization $\beta=j\circ \alpha_X$ where
$j: A\ra A'$ is a morphism of Abelian varieties. The pair $(\alpha, A)$
exists and is unique up to isomorphism for any smooth 
projective variety $X$. The Abelian variety $A$ 
is called the Albanese torus of $X$, denoted by Alb$(X)$.
The universal property implies that $\alpha(X)$ generates
$A$ in the sense that it is not contained in any subabelian
variety and, in particular, $\alpha$ induces an inclusion
$H^0(\Omega_A)\subset H^0(\Omega_X)$.

Let $S_0$ be the minimal model of $S$ and 
$\pi: S\ra S_0$ the projection. Since $\pi$ is a composition of blowups,
its exceptional fibers consist of rational curves that are necessarily
exceptional with respect to $\alpha$. Hence $\alpha=\alpha_0\circ \pi$ 
where $\alpha_0$ is the Albanese map of $S_0$ and Alb$(S)=$Alb$(S_0)$
by the universal property of the Albanese. Since $P^2=K_{S_0}^2\geq 0$ as
 $S_0$ is minimal and since the $\pi$-exceptional curves have zero degree
 with respect to $P=\pi^*K_{S_0}$ and all other irreducible curves in $S$
 are strict transforms of those in $S_0$,  we see that it suffices to prove
 the theorem for $S=S_0$.  So we may set $K=P=(\det d\alpha)$. 
 
 Now let $A$ be an Abelian variety,
 and $\alpha: S\ra A$ a generically finite
 morphism, i.e., $\dim \alpha(S)=2$. 
 Since $S$ has maximal Albanese
 dimension, such a morphism exists.
 Hence,  there is a nonvanishing 
 holomorphic two form $v$ on $A$ 
 such that $\alpha(S)$ is not contained in the codimensional two
 foliation defined by $v$
 and we may set the canonical
 divisor of $S$ to be $K=(\alpha^*v)$. Now $K$ is nef and contains
 all  the exceptional divisors of $\alpha$ by construction. So
 any exceptional divisor of $\alpha$ satisfies the inequality given
 in the theorem. 
 So it remains to prove the inequality of the theorem
 for elliptic curves $C$ that are not $\alpha$-exceptional as
 rational curves are necessarily $\alpha$-exceptional. If $C$
 is such a curve and $\dim A>2$, then $\alpha(C)$ is an elliptic
 curve in $A$ and hence a translate of a one dimensional subgroup
 $E$ of $A$. Let  $\alpha_1$ be the composition of $\alpha$ with the
 projection $p_1$ from $A$ to the quotient abelian variety $A_1=A/E$. If 
 $\dim \alpha_1(S)=2$, then $C$ is $\alpha_1$-exceptional
 and so, as all previous hypothesis on $\alpha$ are satisfied for
 $\alpha_1$, $C$ satisfies the conclusion of our theorem. Hence
 we may assume by induction that $\dim \alpha_1(S)=1$. By 
 the Poincar\'e reducibility theorem, see for example Chap. 6 of \cite{Deb}
 or Theorem 5.3.5 of \cite{BL},
 there is an \'etale base change
 (i.e., an unramified covering map) $z: \tilde A_1\ra A_1$ for an
 Abelian variety such that $z^{-1}(A)=\tilde A_1 \times E$ and
 $z_1^{-1}(p_1): \tilde A_1 \times E \ra A_1$ is the projection 
 $\pi_1$ to the first factor. This 
 means that we have an \'etale covering 
 $z: \tilde A_1 \times E \ra A$ such that 
 $p_1 \circ z=z_1 \circ  \pi_1$. As
 our problem is unchanged by such an
 \'etale base change, we may assume that
 $\tilde A_1=A_1$ and that $A=A_1\times E$ so that $\pi_1$ is the
 quotient map of $A$ by $E$. 
 Let $\tilde C_1$ be the normalization of $C_1=\alpha_1(S)$.
 Then the smooth surface $S'=\tilde C_1\times E$ is the normalization
 of $\alpha(S)$ and, as $S$ is smooth, $\alpha$ factors through $S'$ 
 by the Stein factorization theorem. So replacing $A_1$ by the 
 Albanese torus of $\tilde C_1$, we may assume that $C_1$ is
 smooth.  By construction, $C$ lies in the pre-image of a point 
 $p\in C_1$. Suppose $\dim A_1>1$. Then one can find  a 
 nonzero holomorphic one form $u_1$ on $A_1$ such that
 $T_p C=\ker u_1(p)$. It follows that 
 $\alpha_1^*u_1$ vanishes 
 at every point of $C$ and therefore $C$ lies in the zero locus
 of the section $s$ of ${\cal O}_S(K)=\Omega^2_S$ 
 defined by the wedge product of
 $\alpha_1^*u_1$ with the pull back of a one-form $u_2$ 
 on $E$. Now the pull bak of $u_1$ to $C_1$ is not 
 identically vanishing as $C_1$ generates $A_1$ and so
 $s$ is not identically vanishing. Hence the nef  canonical
 divisor $K=(s)$ contains $C$ giving again
 $$KC\leq K^2\leq 4K^2$$
and so $C$ satisfies the conclusion of the
theorem in this case. We are left with the case 
$C_1=A_1$ so that $\alpha$ is a surjective morphism
to $A$ and so our main theorem reduces to the
following proposition. \\

\begin{prop} With the setup as in proposition~\ref{Pr1}, assume further
that $S$ is a minimal surface. If $C$ is a rational or an elliptic curve in $S$, then
$$KC\leq 4K^2.$$
\end{prop}

 \noindent{\bf Proof:} We will classify a curve as being 
 vertical or horizontal according whether it is $\alpha$-exceptional
 or not and we first note that any element of $F(D)$ for a reduced
 divisor $D$ in $S$ is vertical, being a tree of rational curves
 by lemma~\ref{L1}.
 
  If $C$ is vertical, then $C\leq K$
 and so $KC\leq K^2\leq 4K^2$ by the nefness of $K$. As rational
 curves are necessarily vertical, it suffice to consider the case
 $C$ is a horizontal reduced elliptic curve for our proposition. 
 In this case,  $C_0=\alpha(C)$ is an elliptic curve 
 in $A$. We will also consider $C_0$ as a reduced divisor in $A$. 
 Let $C'$ be the sum of the non-elliptic  
 horizontal components of $\hat C=\alpha^*C_0$ 
 on which $\hat C$ has multiplicity one and let $\bar C = \alpha^*C_0 -C'$. 
 Then $C \leq \bar C\leq  \alpha^*C_0$ and so 
 $$KC\leq K\bar C$$ as $K$ is nef. 
 Let $\sum_i n_i D_i$ be the prime decomposition of $K$.  
 We are reduced to bounding the intersection number 
 $$K\bar C=\sum_i  n_i(D_i, \bar C).$$

 Let $D_i$ be a component of $K$ that is either vertical or $D_i\leq \bar C$ is 
 horizontal,  then $$(D_i, \bar C)\leq (D_i, \alpha^*C_0)=0$$ 
 by the definition of $\bar C$ 
 (since $C'$ is an effective horizontal divisor not containing $D_i$).\\

 We now set $D=K_{red}$ in what to follow and decompose $D=\sum_i D_i$  into three parts:
 $$D=D^0+D'+D'',\ \ \mbox{where}$$
 
 \begin{itemize}
 
 \item $D^0$ consists of all the vertical components of $D$ connected
 via vertical curves 
 to a part of $D$ that belongs to $F(D)\setminus F^{HZ}(D)$ or to a part
 that is a maximal element of $F^{HJ}(D)$ not attached to any horizontal
 component of $D$,
 
 \item $D'$ consists of horizontal components of $D$ not in $\bar C$
 and of any element of $F^{HJ}(D)$ that is attached to them,
 
 \item $D''$ consists of horizontal components of $D$ in $\bar C$
 and of any element of $F^{HJ}(D)$ that is attached to them.
 
 \end{itemize}

 We first observe that there are no maximal element of $F(D)$ that is isolated
 as that would lead to an isolated surface singularity $p$ 
 sitting above a smooth surface but away from any ramification 
 divisor, an impossibility by theorem III.5.2 of \cite{BPV}. We observe
 also, by our last inequality above, that $K_{D^0}$ and $K_{D''}$ have
 non-positive intersection with $\bar C$. \\

 By definition, we can write $D^0=\sum_i C^M_i+\sum_j (C^{HJ}_j+C^+_j)+C'$ 
 where the $C^M_i$'s are the maximal elements in $F(D)$ not in $F^{HJ}(D)$, 
 the $C^{HJ}_j$'s are the maximal elements in $F^{HJ}(D)$ not attached to any 
 horizontal part of $D$, $C^+_j$ the vertical component of $D$ attached 
 to $C^{HJ}_j$ and $C'$ the rest of $D^0$. By construction, $C'$ does 
 not contain any leave of $F_D$ so that any component $D_k$ of $C'$ 
 satisfies $(D_k, K+D)\geq 0$. Hence $(K_{C'}, K+D)\geq 0$, where 
 we recall that $K_{C'}$ is the part of $K$ supported on $C'$. 
 We now note that, for all $j$,  $\hat C_j=C^{HJ}_j+C^+_j\not\in F(D)$ so that 
 $(\hat C_j,K+D)\geq 0$. Letting $m_j$ be the multiplicity of $K$ on the leave 
 of $C^{HJ}_j$, we have $m_j\hat C_j\leq K$ by lemma~\ref{L2}. As 
 $K_{\hat C_j}-m_j\hat C_j$ consists of components not belonging to 
 $F(D)$, we obtain $(K_{\hat C_j}, K+D)\geq 0$ for all $j$. Finally, as 
 $(K_{C^M_i}, K+D)\geq 0$ by lemma~\ref{L3} and as
 $K_{D^0}=\sum_i K_{C^M_i} + \sum_j K_{\hat C_j} +K_{C'}$, we get 
 $$(K_{D^0}, \bar C)\leq 0\leq (K_{D^0}. K+D)\leq 2(K_{D^0}. K+D).$$

 By the same token, $D''\not\in F(D)$. In fact, by definition,
 $D''=\sum_l (D_l+\sum_{l'} D_{ll'})$ where the $D_l$'s are
 the horizontal components of $D$ lying in $\bar C$ and, for
 each $l$, the $D_{ll'}$'s are the elements of $F^{HJ}(D)$ attached to $D_l$
 and $E_l=D_l+\sum_{l'} D_{ll'}\not\in F(D)$. Let $m_l$ be the multiplicity
 of $K$ on $D_l$. Then $m_l$ is greater than the multiplicity $m_{ll'}$ of $K$
 on the leave of $D_{ll'}$  for every $l'$ by lemma~\ref{L2}. We also have 
 $\hat D_l= m_l D_l +\sum_{l'} m_{ll'}D_{ll'}\leq K_{D''}$
 by the same lemma. As $(m_lE_l, K+D)\geq 0$ while $(D_{ll'}, K+D)<0$ for all $l'$,
 we see that $(\hat D_l, K+D)\geq 0$ for all $l$. Since 
 $K_{D''}-\sum_l \hat D_l$ is an effective divisor not having
 any component that belongs to $F(D)$,  we obtain as before 
 $$(K_{D''}, \bar C)\leq 0\leq (K_{D''}, K+D)\leq 2(K_{D''}, K+D).$$
 
 As for $D'$,  let 
 $\{G_n\}_{n\in I}$ be the collection of 
 horizontal components of $D$ not lying in $\bar C$ and,
 for each $n\in I$,
 let  $H_n=G_n+G_{nn'}$  where the $G_{nn'}$'s are the elements of 
 $F^{HJ}(D)$ attached to $G_n$.
 Then $D'=\sum_n H_n$ and we are interested in bounding
 $(K_{H_n},\bar C)$ for each $n$. Fix an $n$.
 We have by proprosition~\ref{Pr1} that 
 $$ (G_n,\bar C)\leq 2(G_n, K+D) - n_{G_n},$$
 where $n_{G_n}$ is the intersection number of 
 the smooth part of $G$
 with the rest of $D$. 
 Since $(G_{nn'}, G_n)=1=-(G_{nn'}, K+D)$, we have
 $n_{G_n}\geq (\sum_{n'} G_{nn'}, G_n)=-(\sum_{n'} G_{nn'}, K+D)$
 and so
 $$ (G_n,\bar C)\leq 2(G_n, K+D) - n_{G_n}\leq 2(G_n, K+D) + (\sum_{n'} G_{nn'}, K+D). $$
 Let $m_n$ be the multiplicity of $K$ on $G_n$ and $m_{nn'}$ that on the leave of
 $G_{nn'}$. Then lemma~\ref{L2} implies that $2m_{nn'}\leq m_n$ for all $n'$ and
 that $\bar H_n=m_nG_n + \sum_{n'} m_{nn'}G_{nn'}\leq K_{H_n}$. It follows from
 $(G_{nn'}, \bar C)\leq 0$ and $(G_{nn'}, K+D)<0$ that 
 $$(\bar H_n,\bar C)\leq (m_nG_n,\bar C)\leq 
 (2m_nG_n + \sum_{n'} m_nG_{nn'}, K+D)
 \leq (2\bar H_n, K+D)=2(\bar H_n, K+D). $$
 Since any component $D_r$ of  $K_{H_n}-\bar H_n$ is a vertical non-leave
 component of $D$, we have $(D_r, \bar C)\leq 0\leq 2(D_r, K+D)$ and so
 $ (K_{H_n}, \bar C)\leq 2(K_{H_n}, K+D).$ Summing over $n$ gives
 $$ (K_{D'},\bar C)\leq 2 (K_{D'},K+D).$$
 
 Finally, as $K=K_{D^0}+K_{D'}+K_{D''}$, we obtain
 $$(K,\bar C)\leq 2 (K,K+D)\leq 4K^2.\ \ \BOX$$


\medskip
\noindent
Steven Shin-Yi Lu\\
D\'epartment de Math\'ematiques\\
UQAM\\
C.P. 8888 Succursale Centre-ville\\
Montr\'eal   H3C 3P8\\
Canada\\

\noindent
lu.steven@uqam.ca\\


\begin{thebibliography}{99}

 \bibitem{Ab} D. Abramovich, 
Subvarieties of semiabelian varieties.  
Compositio Math.  {\bf 90}  (1994),  no. 1, 37--52. 

\bibitem{BL} C. Birkenhake,  H. Lange, 
Complex abelian varieties. 
Second edition. Grundlehren der Mathematischen Wissenschaften 
[Fundamental Principles of Mathematical Sciences], 
{\bf 302}. Springer-Verlag, Berlin, 2004. 635 pp.

 \bibitem{BPV} W. Barth, C. Peters,  A. Van de Ven,  
Compact complex surfaces. 
Ergebnisse der Mathematik und ihrer Grenzgebiete (3) 
[Results in Mathematics and Related Areas (3)], {\bf 4}. 
Springer-Verlag, Berlin, 1984. 

 \bibitem{Bog} F. A. Bogomolov,  
Holomorphic tensors and vector bundles on projective manifolds. (Russian)  
Izv. Akad. Nauk SSSR Ser. Mat.  {\bf 42}  (1978), no. 6, 1227--1287, 1439. 

 \bibitem{Deb} O. Debarre, Tores et vari\'et\'s ab\'eliennes complexes. 
 [Complex tori and abelian varieties] Cours Sp\'ecialis\'es, {\bf 6}. 
 Soci\'et\'e Math\'ematique de France, Paris; EDP Sciences, Les Ulis, 1999. 125 pp.

 \bibitem{GG} M. Green, P. Griffiths, 
Two applications of algebraic geometry to entire holomorphic mappings.  
The Chern Symposium 1979 (Proc. Internat. Sympos., Berkeley, Calif., 1979),  
pp. 41--74, Springer, New York-Berlin, 1980. 

 \bibitem{Kaw} Y. Kawamata,
Characterization of abelian varieties.  
Compositio Math.  {\bf 43}  (1981), no. 2, 253--276. 

 \bibitem{Kaw1} Y. Kawamata,  
On Bloch's conjecture.  Invent. Math.  {\bf 57}  (1980), no. 1, 97--100.

 \bibitem{KM} J. Koll\`ar, T. Matsusaka,
Riemann-Roch type inequalities.
Amer. J. Math. {\bf 105} (1983), no. 1, 229--252. 

 \bibitem{Lang}
S. Lang, Hyperbolic and Diophantine analysis.  
Bull. Amer. Math. Soc. (N.S.)  {\bf 14}  (1986),  no. 2, 159--205.

 \bibitem{Langer}  A. Langer,  
Logarithmic orbifold Euler numbers of surfaces with applications. 
Proc. London Math. Soc. (3)  {\bf 86}  (2003),  no. 2, 358--396. 

 \bibitem{LM1}
S. Lu, Y. Miyaoka, 
Bounding curves in algebraic surfaces by genus and Chern numbers.  
Math. Res. Lett.  {\bf 2}  (1995),  no. 6, 663--676.

 \bibitem{LM2}
S. Lu, Y. Miyaoka, 
Bounding codimension-one subvarieties 
and a general inequality between Chern numbers.  
Amer. J. Math.  {\bf 119}  (1997),  no. 3, 487--502. 

 \bibitem{Lu10} S. Lu,  On varieties of general type
with maximal Albanese dimension,  Preprint.

 \bibitem{Miya} Y. Miyaoka, 
On the Chern numbers of surfaces of general type.  
Invent. Math.  {\bf 42}  (1977), 225--237. 

 \bibitem{Miya07} Y. Miyaoka, The orbifold Miyaoka-Yau
 Inequality and an effective Bogomolov-McQuillan Theorem. Preprint 2007.

 \bibitem{Nog86}
J. Noguchi,  
Logarithmic jet spaces and extensions of de Franchis' theorem.  
Contributions to several complex variables,  
227--249, Aspects Math., E9, Vieweg, Braunschweig, 1986.

 \bibitem{NWY}
J. Noguchi, J. Winkelmann, K. Yamanoi,  
The second main theorem for holomorphic curves into semi-abelian varieties.  
Acta Math.  {\bf 188}  (2002),  no. 1, 129--161.

 \bibitem{NWY07}
J. Noguchi, J. Winkelmann, K. Yamanoi, 
Degeneracy of holomorphic curves into algebraic varieties.  
J. Math. Pures Appl. (9)  {\bf 88}  (2007),  no. 3, 293--306. 

 \bibitem{Sakai} F.  Sakai,  
Semistable curves on algebraic surfaces 
and logarithmic pluricanonical maps.  
Math. Ann.  {\bf 254}  (1980), no. 2, 89--120. 

 \bibitem{Ueno} K. Ueno, Classification theory of algebraic varieties and compact complex spaces.  Lecture Notes in Mathematics, {\bf 439}. Springer-Verlag, Berlin-New York, 1975. 278 pp.

 \bibitem{Yau} S.-T. Yau, 
On the Ricci curvature of a compact 
K\"ahler manifold and the complex Monge-Amp\`ere equation. I.  
Comm. Pure Appl. Math.  {\bf 31}  (1978), no. 3, 339--411. 

 \bibitem{Zariski}
O. Zariski, 
The theorem of Riemann-Roch for high multiples 
of an effective divisor on an algebraic surface.  
Ann. of Math. (2)  {\bf 76}  (1962) 560--615. 

\end{thebibliography}
\end{document}